# Analytic perturbations and systematic bias in statistical modeling and inference[*]

## Jerzy A. Filar[1], Irene Hudson[1], Thomas Mathew[2] and Bimal Sinha[2]

*University of South Australia and University of Maryland Baltimore County*

**Abstract:** In this paper we provide a comprehensive study of statistical inference in linear and allied models which exhibit some analytic perturbations in their design and covariance matrices. We also indicate a few potential applications. In the theory of perturbations of linear operators it has been known for a long time that the so-called "singular perturbations" can have a big impact on solutions of equations involving these operators even when their size is small. It appears that so far the question of whether such undesirable phenomena can also occur in statistical models and their solutions has not been formally studied. The models considered in this article arise in the context of nonlinear models where a single parameter accounts for the nonlinearity.

## 1. Introduction

The problem of estimating parameters and drawing suitable inference about them from noisy data that follow classical linear models has a long history in experimental science and engineering. A standard assumption in such models is the normality and independence of the random error terms, and a complete and certain knowledge of the associated design matrices. Since there may be situations where such assumptions are not valid, or satisfactory, it is of interest to develop appropriate modifications of the usual standard solutions. In this regard, some robustness studies have been reported in the literature. Kariya and Sinha [12] discuss distributional robustness in terms of deviations from normality of the error terms. They show that, under quite general conditions, the standard inference for the fixed effects under normality based on the $F$ test continues to hold under a broad class of elliptically symmetric distributions. Calafiore and Ghaout [5] discuss some aspects of linear models in the presence of uncertainties in the mean and covariance matrix of the observed data, and study robust estimation of the underlying parameters. We also refer to Chandrasekaran et al. [6], Ghaoui and Lebret [9] for some related work.

Our goal in this paper is to study some aspects of analytic perturbations in statistical modeling and inference in the context of linear and allied models in the presence of certain "structured and systematic uncertainties" in the design and covariance matrices. With respect to these uncertainties in the design matrices in linear models, our approach is different from the usual errors-in-variables models in that the uncertainties in the design matrices in our context are non-stochastic in nature and may be a result of small (but unknown) measurement errors that persistently underestimate or overestimate the true values of certain variables. Note that these kinds of deviations from true measurements are not random as they may be caused either by a persistent deficiency in the measurement sensors or even by consistent attempts at deception in a small (and hence hard to detect?) way. We propose the name *systematic measurement bias* to describe this class of uncertainties in the design and covariance matrices. The scenario considered here may also arise when we have nonlinear models where a single parameter accounts for the nonlinearity; in other words, the model is linear if this particular parameter is zero.

---

[*]Supported, in part, by Australian Research Council Grants DP0343O28 and DP0774504.
[1]School of Mathematics and Statistics, University of South Australia, Mawson Lakes Campus, Mawson Lakes, SA 5095, Australia, e-mail: Jerzy.Filar@unisa.edu.au; Irene.Hudson@unisa.edu.au
[2]Department of Mathematics and Statistics, University of Maryland Baltimore County, Baltimore, Maryland 21250, USA, e-mail: mathew@math.umbc.edu; sinha@math.umbc.edu

Our ability to identify, quantify and analyze this kind of a bias stems from certain recent developments in the well established subject of *analytic perturbation theory* of matrices and operators. The reader is referred to the famous treatise by Kato [13] for a comprehensive treatment of this subject. However, in this study, results obtained by Avrachenkov [2], Avrachenkov and Haviv [3] and Avrachenkov et al. [4] (see Section 2) provide us with tools to assess the impact of these systematic measurement biases on the quality of the estimates and the inferences drawn from these models. Note that in the theory of perturbations of operators it has been known for a long time that some small perturbations can have a big impact on solutions of equations involving these operators; these are the so-called "singular perturbations". The question of whether such undesirable phenomena can also occur in statistical models and their solutions appears not to have been formally studied so far. Our application in the area of principal component analysis (PCA) involves perturbations in the estimated sample covariance matrix via perturbations in the underlying design matrices while the application in a factor analysis model deals with perturbations in the covariance matrix.

The organization of the paper is as follows. We introduce some background of perturbation theory and give associated Laurent series and fundamental equations in Section 2. The context of linear models with perturbations in the design matrices is developed in Section 3, a core section of the paper. We derive the necessary expansions of the associated estimates of the parameters and their covariance matrices, and provide an outline of the inference on estimable linear parametric functions under perturbations of the underlying design matrices. Perturbations in principal component analysis appear in Section 4, and Section 5 deals with perturbations in the context of factor analysis models. Some applications are indicated in Section 6.

## 2. Background: Laurent series and fundamental equations

In this section we briefly review some results concerning analytic perturbations of matrices. We refer the reader to Avrachenkov [2] and Avrachenkov et al. [4] for details of proofs and derivations of these results.

Let $\{A_k\}_{k=0,1,\ldots} \subseteq C^{n \times n}$ be a sequence of matrices that defines an analytic matrix valued function

$$(2.1) \qquad A(\varepsilon) = A_0 + \varepsilon A_1 + \varepsilon^2 A_2 + \cdots.$$

The above series is assumed to converge in some non-empty neighbourhood of $\varepsilon = 0$. We say that $A(\varepsilon)$ is an *analytic perturbation* of the matrix $A_0 = A(0)$, and that it is a *regular perturbation* if $A_0$ is nonsingular, and a *singular perturbation* otherwise. Even though all results of this paper can be extended to the above general analytic perturbation, the statistical considerations are sufficiently well illustrated by the somewhat simpler case of a *linear perturbation* where (2.1) is replaced by

$$(2.2) \qquad A(\varepsilon) = A + \varepsilon B,$$

where $A_0 = A$, $A_1 = B$, and $A_k = 0$, $k \geq 2$. Assume the inverse matrices $A^{-1}(\varepsilon)$ exist in some (possibly punctured) disc centred at $\varepsilon = 0$. Beyond linear perturbations, we may also have a quadratic or a polynomial perturbation though the contributions of the higher order terms will in general be small. Let us consider the following example of an analytically (linear) perturbed singular matrix

$$A(\varepsilon) = \begin{bmatrix} 1-\varepsilon & 1+\varepsilon \\ 1-2\varepsilon & 1-\varepsilon \end{bmatrix} = \begin{bmatrix} 1 & 1 \\ 1 & 1 \end{bmatrix} + \varepsilon \begin{bmatrix} -1 & 1 \\ -2 & -1 \end{bmatrix}.$$

Its inverse is given by

$$A^{-1}(\varepsilon) = \frac{1}{-\varepsilon(1-3\varepsilon)} \begin{bmatrix} 1-\varepsilon & -1-\varepsilon \\ -1+2\varepsilon & 1-\varepsilon \end{bmatrix}.$$

It should be clear that the above inverse admits a Laurent series expansion. To see this, we just expand $\det A(\varepsilon))^{-1} = 1/(-\varepsilon(1-3\varepsilon))$ as a scalar power series in $\varepsilon$, multiply it by $\operatorname{adj} A(\varepsilon)$ and equate coefficients with the same power of $\varepsilon$. In this case, we have

$$\begin{aligned} A^{-1}(\varepsilon) &= (-\frac{1}{\varepsilon} - 3 - 9\varepsilon - \cdots) \begin{bmatrix} 1-\varepsilon & -1-\varepsilon \\ -1+2\varepsilon & 1-\varepsilon \end{bmatrix} \\ &= \frac{1}{\varepsilon}\begin{bmatrix} -1 & 1 \\ 1 & -1 \end{bmatrix} + \begin{bmatrix} -2 & 4 \\ 1 & -2 \end{bmatrix} + \varepsilon \begin{bmatrix} -6 & 12 \\ 3 & -6 \end{bmatrix} + \cdots. \end{aligned}$$

which is obviously a Laurent series. Clearly, the above is a very inefficient method of deriving this Laurent series. Fortunately, in the case of a linear perturbation, the following result that follows from the analysis in Avrachenkov [2] and Avrachenkov et al. [4] provides a much better method.

If $A^{-1}(\varepsilon)$ exists in a punctured neighborhood of $\varepsilon = 0$, then it can be expanded as a Laurent series

$$(2.3) \qquad A^{-1}(z) = A_S^{-1}(\varepsilon) + A_R^{-1}(\varepsilon) = \frac{1}{\varepsilon^s} Y_{-s} + \cdots + \frac{1}{\varepsilon} Y_{-1} + Y_0 + \varepsilon Y_1 + \cdots,$$

where $A_S^{-1}(\varepsilon)$ is the *singular component* of $A^{-1}(\varepsilon)$ and consists of the part of the series made of terms involving negative powers of $\varepsilon$, and $A_R^{-1}(\varepsilon)$ is the *regular component* consisting of all the remaining terms. The coefficients $Y_k$, $k = -s, -s+1, \ldots$ satisfy the following recursions

$$(2.4) \qquad Y_{k+1} = (-Y_0 B) Y_k, \quad k = 0, 1, \ldots,$$

$$(2.5) \qquad Y_{-k-1} = (-Y_{-1} A) Y_{-k}, \quad k = 1, \ldots, s-1.$$

Thus the knowledge of the matrices $A, B, Y_0$ and $Y_{-1}$ is sufficient to determine all coefficients of the Laurent series expansion (2.3).

There are now two cases to consider. In the case of a regular linear perturbation, $A_S^{-1}(\varepsilon) = 0 = Y_{-1} = Y_{-k}$, $k = 2, \ldots s$, and $Y_0 = A^{-1} = A_0^{-1}$. However the singularly perturbed case is a little more complicated even when the perturbation is linear. Firstly, the order of the singularity, the index $s$ in (2.3), which is also known as the *pole* of the expansion of $A^{-1}(\varepsilon)$, needs to be determined. Secondly, the identity $A(\varepsilon)A^{-1}(\varepsilon) = I$ provides a set of *fundamental equations* obtained by equating the coefficients of like powers of $\varepsilon$ that supplies the unique set of coefficients of the expansion (2.3). Fortunately, only $s + 1$ of these are needed to determine $Y_{-s}, Y_{-s+1} \cdots Y_{-1}, Y_0$ uniquely.

In particular, for each $t = 0, 1, \ldots$ let us define the following augmented matrix

$$\mathcal{A}^{(t)} = \begin{bmatrix} A_0 & 0 & 0 & \cdots & 0 \\ A_1 & A_0 & 0 & \cdots & 0 \\ A_2 & A_1 & A_0 & \cdots & 0 \\ \vdots & \vdots & \vdots & \ddots & \vdots \\ A_t & A_{t-1} & \cdots & A_1 & A_0 \end{bmatrix}.$$

Of course, in the case of a linear perturbation the above augmented matrices have the form

$$\mathcal{A}_L^{(t)} = \begin{bmatrix} A & 0 & 0 & \cdots & 0 \\ B & A & 0 & \cdots & 0 \\ 0 & B & A & \cdots & 0 \\ \vdots & \vdots & \vdots & \ddots & \vdots \\ 0 & 0 & \cdots & B & A \end{bmatrix}.$$

The determination of the order of the pole $s$ in (2.3) can be achieved with the help of the following result that is proved in Avrachenkov [2] and Avrachenkov et al. [4].

**Theorem 2.1.** *The order of the pole $s$ is given by the smallest value of $t$ for which* $\operatorname{rank}[\mathcal{A}^{(t)}] = \operatorname{rank}[\mathcal{A}^{(t-1)}] + n$, *where $n$ is the dimension of $A(\varepsilon)$.*

Now if we define two block *column matrices* $\mathcal{Y} := [Y_{-s}, \ldots, Y_{-1}, Y_0]^T$ and $\mathcal{J} := [0, \ldots, 0, I]^T$ and let $s$ be the order of the pole of the expansion of $A^{-1}(\varepsilon)$, then the coefficients of the non-positive powers of $\varepsilon$ in (2.3) are given as the unique solution of the linear system:

$$(2.6) \qquad \mathcal{A}^{(s)} \mathcal{Y} = \mathcal{J}.$$

Next, we obtain a recursive formula for the Laurent series coefficients in (2.3). To eliminate negative indices in subscripts, it will be convenient to define $X_k := Y_{k-s}$, $k \geq 0$ and to introduce the Moore-Penrose generalized inverse of the augmented matrix $\mathcal{A}^{(s)}$. In particular, let $\mathcal{G}^{(s)} \stackrel{def}{=} [\mathcal{A}^{(s)}]^+$ be the Moore-Penrose generalized inverse of $\mathcal{A}^{(s)}$ and define the submatrices $G_{ij}^{(s)} \in C^{n \times n}$ for $0 \leq i, j \leq t$ by

$$\mathcal{G}^{(s)} = \begin{bmatrix} G_{00}^{(s)} & \cdots & G_{0s}^{(s)} \\ \vdots & \ddots & \vdots \\ G_{s0}^{(s)} & \cdots & G_{ss}^{(s)} \end{bmatrix},$$

where the dimensions and locations of $G_{ij}^{(s)}$ are in correspondence with the block structure of $\mathcal{A}^{(s)}$.

In Avrachenkov et al. [4] it is shown that the first $n$ rows of the generalized inverse $\mathcal{G}^{(s)}$, namely, $[G_{00}^{(s)} \cdots G_{0s}^{(s)}]$ are all that is needed to calculate the coefficients of the expansion (2.3). Indeed, the following recursive formula provides the solution:

$$(2.7) \qquad X_k = \sum_{j=0}^{s} G_{0j}^{(s)}(\delta_{j+k,s}I - \sum_{i=1}^{k} A_{i+j}X_{k-i}), \quad k = 1, 2, \ldots,$$

initialising with $X_0 = G_{0s}^{(s)}$.

To illustrate these formulae we next use them to verify the first two coefficients in the expansion of $A^{-1}(\varepsilon)$. First, we note that $\mathcal{A}_L^{(0)} = A$, and

$$\mathcal{A}_L^{(1)} = \begin{bmatrix} A & 0 \\ B & A \end{bmatrix}.$$

It is now easy to check that $\text{rank}[\mathcal{A}^{(1)}] = 3 = \text{rank}[\mathcal{A}^{(0)}] + n = 1 + 2$. Thus, by Theorem 2.1, the order of the pole is indeed $s = 1$ and

$$\mathcal{G}^{(1)} = \begin{bmatrix} 1 & 1 & -1 & 1 \\ -.5 & -.5 & 1 & -1 \\ .75 & .75 & -.5 & 1 \\ .75 & .75 & -.5 & 1 \end{bmatrix}.$$

It immediately follows from the notation introduced earlier and (2.7) that

$$Y_{-1} = X_0 = G_{01}^{(1)} = \begin{bmatrix} -1 & 1 \\ 1 & -1 \end{bmatrix}.$$

Now, using (2.7) again and the fact that $A_k = 0$ for $k \geq 2$, we obtain

$$Y_0 = X_1 = G_{00}^{(1)}[I_2 - BX_0] + G_{01}^{(1)}[0 - 0X_0],$$

that is,

$$Y_0 = \begin{bmatrix} 1 & 1 \\ -.5 & -.5 \end{bmatrix} \begin{bmatrix} \begin{bmatrix} 1 & 0 \\ 0 & 1 \end{bmatrix} - \begin{bmatrix} -1 & 1 \\ -2 & -1 \end{bmatrix} \begin{bmatrix} -1 & 1 \\ 1 & -1 \end{bmatrix} \end{bmatrix} = \begin{bmatrix} -2 & 4 \\ 1 & -2 \end{bmatrix}.$$

It is, perhaps, worth commenting that while, at first sight, the dimension of $\mathcal{A}^{(s)}$ may appear prohibitively large, it appears that the singularity of order $s = 1$ occurs very frequently (this is made precise in Avrachenkov [2]).

## 3. Perturbed linear model

We begin with the standard linear model

$$(3.1) \qquad \mathbf{Y}_{n \times 1} = X^T(\varepsilon)_{n \times m} \boldsymbol{\beta}_{m \times 1} + \boldsymbol{\xi}_{n \times 1}$$

and assume that the underlying design matrix $X(\varepsilon)$ admits an analytic perturbation expansion of the form

$$(3.2) \qquad X(\varepsilon) = X_0 + \varepsilon X_1 + \varepsilon^2 X_2 + \varepsilon^3 X_3 + \cdots$$

22

where $X_0, X_1, X_2, \ldots$ are known *component* design matrices, and the error terms $\boldsymbol{\xi} \sim N[0, \sigma^2]$, which is the usual assumption in linear models. Here $\varepsilon > 0$ is referred to as a *perturbation* parameter which is typically small and also unknown. As mentioned above, our goal in this section is to study the effects of the perturbation parameter $\varepsilon$ on the statistical inference about $\boldsymbol{\beta}$. In view of the representation of $X(\varepsilon)$, the relevance of the discussion in Section 2 is obvious.

Towards this end, we first note that for a given $\epsilon$, the maximum likelihood (also the least squares) estimate of $\boldsymbol{\beta}$, which is obtained by minimizing $[\mathbf{Y} - X^T(\varepsilon)\boldsymbol{\beta}]'[\mathbf{Y} - X^T(\varepsilon)\boldsymbol{\beta}]$ with respect to $\boldsymbol{\beta}$, is given by

$$(3.3) \qquad \hat{\boldsymbol{\beta}}(\varepsilon) = [X(\varepsilon)X^T(\varepsilon)]^{-1}X(\varepsilon)\mathbf{Y}.$$

Moreover, an estimate of the error variance $\sigma^2$ is obtained from the residual sum of squares $SSE(\varepsilon)$ given by

$$(3.4) \qquad \begin{aligned} SSE(\varepsilon) &= \mathbf{Y}^T[I_n - X^T(\varepsilon)\{X(\varepsilon)X^T(\varepsilon)\}^{-1}X(\varepsilon)\mathbf{Y} \\ &= \mathbf{Y}^T P(\varepsilon) \mathbf{Y} \end{aligned}$$

where $P(\varepsilon)$, the usual projection operator, is given by

$$(3.5) \qquad P(\varepsilon) = [I_n - X^T(\varepsilon)\{X(\varepsilon)X^T(\varepsilon)\}^{-1}X(\varepsilon)].$$

To study the dependence of $\hat{\boldsymbol{\beta}}(\varepsilon)$ and $SSE(\varepsilon)$ on $\varepsilon$, let us recall from (3.2) that $X(\varepsilon) = X_0 + \varepsilon X_1 + \varepsilon^2 X_2 + \cdots$, which yields

$$(3.6) \qquad \begin{aligned} X(\varepsilon)X^T(\varepsilon) &= X_0 X_0^T + \varepsilon\{X_0 X_1^T + X_1 X_0^T\} \\ &\quad + \varepsilon^2\{X_0 X_2^T + X_1 X_1^T + X_2 X_0^T\} \\ &\quad + \cdots \\ &= B_0 + \varepsilon B_1 + \varepsilon^2 B_2 + \cdots \end{aligned}$$

where $B_0 = X_0 X_0^T$ and $B_i$ is symmetric, $\forall i$. We now distinguish between two cases depending on whether $B_0$ is nonsingular or singular.

Case 1: $B_0$ is nonsingular. In this case which corresponds to a regular perturbation of $X(\varepsilon)$, we readily get for small $\varepsilon$

$$(3.7) \qquad \begin{aligned} (X(\varepsilon)X^T(\varepsilon))^{-1} &= (B_0 + \varepsilon B_1 + \varepsilon^2 B_2 + \cdots)^{-1} \\ &= C_0 + \varepsilon C_1 + \varepsilon^2 C_2 + \cdots \end{aligned}$$

where $C_0 = B_0^{-1} \leftrightarrow C_0^{-1} = B_0$ and the remaining $C_i$'s can be computed following the ideas of Section 2. It then follows from (3.3) that

$$(3.8) \qquad \begin{aligned} \hat{\boldsymbol{\beta}}(\varepsilon) &= (C_0 + \varepsilon C_1 + \varepsilon^2 C_2 + \cdots)(X_0 + \varepsilon X_1 + \cdots)\mathbf{Y} \\ &= C_0 X_0 \mathbf{Y} + \varepsilon[C_1 X_0 + C_0 X_1]\mathbf{Y} \\ &\quad + \varepsilon^2[C_2 X_0 + C_1 X_1 + C_0 X_2]\mathbf{Y} + \cdots \\ &= \hat{\boldsymbol{\beta}} + \varepsilon[C_1 C_0^{-1}\hat{\boldsymbol{\beta}} + C_0 X_1 \mathbf{Y}] \\ &\quad + \varepsilon^2[C_2 C_0^{-1}\hat{\boldsymbol{\beta}} + C_1 1 X_1 \mathbf{Y} + C_0 X_2 \mathbf{Y}] \\ &\quad + \cdots \\ &= [C_0 + \varepsilon C_1 + \varepsilon^2 C_2 + \cdots]C_0^{-1}\hat{\boldsymbol{\beta}} \\ &\quad + \varepsilon[C_0 X_1 \mathbf{Y}] + \varepsilon^2[C_1 X_1 \mathbf{Y} + C_0 X_2 \mathbf{Y}] + \cdots \\ &\approx \hat{\boldsymbol{\beta}} + \varepsilon[C_1 B_0 \hat{\boldsymbol{\beta}} + C_0 X_1 \mathbf{Y}] \text{ for small } \varepsilon \end{aligned}$$

where

$$\hat{\boldsymbol{\beta}} = B_0^{-1} X_0 \mathbf{Y} = C_0 X_0 \mathbf{Y}$$

is the usual regression coefficient estimate without any perturbation parameter. Moreover, from (3.7), we get

$$
\begin{aligned}
(3.9) \quad & X^T(\varepsilon)(X(\varepsilon)X^T(\varepsilon))^{-1}X(\varepsilon) \\
&= (X_0 + \varepsilon X_1 + \varepsilon^2 X_2 + \cdots)^T \\
&\quad \times (C_0 + \varepsilon C_1 + \varepsilon^2 C_2 + \cdots)(X_0 + \varepsilon X_1 + \cdots) \\
&= X_0^T C_0 X_0 + \varepsilon[X_1^T C_0 X_0 + X_0^T C_1 X_0 + X_0^T C_0 X_1] \\
&\quad + \varepsilon^2 [X_2^T C_0 X_0 + X_0^T C_2 X_0 + X_0^T C_0 X_2 + X_1^T C_1 X_0 \\
&\quad\quad + X_0^T C_1 X_1 + x_1^T C_0 X_1] \\
&\quad + \cdots \\
&\approx X_0^T C_0 X_0 + \varepsilon[X_1^T C_0 X_0 + X_0^T C_1 X_0 + X_0^T C_0 X_1] \text{ for small } \varepsilon
\end{aligned}
$$

which readily yields the following expansion of the projection operator

$$
\begin{aligned}
(3.10) \quad P(\varepsilon) &= I_n - X^T(\varepsilon)\{X(\varepsilon)X^T(\varepsilon)\}^{-1}X(\varepsilon) \\
&= [I_n - X_0^T C_0 X_0] - \varepsilon[X_1^T C_0 X_0 + X_0^T C_1 X_0 + X_0^T C_0 X_1] \\
&\quad - \varepsilon^2 [X_2^T C_0 X_0 + X_0^T C_2 X_0 + X_0^T C_0 X_2 + X_1^T C_1 X_0 \\
&\quad\quad + X_0^T C_1 X_1 + X_1^T C_0 X_1] \\
&\quad - \cdots
\end{aligned}
$$

so that, from (3.4), we get

$$
\begin{aligned}
(3.11) \quad SSE(\varepsilon) &= SSE - \varepsilon \mathbf{Y}^T [X_1^T C_0 X_0 + X_0^T C_1 X_0 + X_0^T C_0 X_1] \mathbf{Y} \\
&\quad - \varepsilon^2 \mathbf{Y}^T \left[X_2^T C_0 X_0 + X_0^T C_2 X_0 + X_0^T C_0 X_2 + X_1^T C_1 X_0 \right. \\
&\quad\quad \left. + X_0^T C_1 X_1 + X_1^T C_0 X_1\right] \mathbf{Y} - \cdots
\end{aligned}
$$

where $SSE = \mathbf{Y}^T(I_n - X_0^T C_0 X_0)\mathbf{Y}$ is the standard error sum of squares without any perturbation parameter. The expansions in (3.8) and (3.11) clearly reveal the effects of the perturbation parameter $\varepsilon$ on the estimates of $\boldsymbol{\beta}$ and $\sigma^2$.

We now turn our attention to the problem of testing $H_0 : \boldsymbol{\beta} = \boldsymbol{\beta}_0$ vs $H_1 : \boldsymbol{\beta} \neq \boldsymbol{\beta}_0$. Under normality and independence of the error terms, a standard test in this context is the familiar $F$-test based on the $F$-statistic given by

$$
\begin{aligned}
F(\varepsilon) &= \frac{(\hat{\boldsymbol{\beta}}(\varepsilon) - \boldsymbol{\beta}_0)^T [X(\varepsilon)X^T(\varepsilon)](\hat{\boldsymbol{\beta}}(\varepsilon) - \boldsymbol{\beta}_0)/m}{SSE(\varepsilon)/(n-m)} \\
(3.12) \quad &\sim F_{m,n-m}(H_0).
\end{aligned}
$$

To study the dependence of $F(\varepsilon)$ on $\varepsilon$, we proceed as follows. From (3.8) and (3.11),

we get

$$
\begin{aligned}
(3.13) \quad & (\hat{\boldsymbol{\beta}}(\varepsilon)) - \boldsymbol{\beta}_0)^T [X(\varepsilon) X^T(\varepsilon)] (\hat{\boldsymbol{\beta}}(\varepsilon) - \boldsymbol{\beta}_0) \\
&= [(\hat{\boldsymbol{\beta}} - \boldsymbol{\beta}_0) + \varepsilon (C_1 X_0 + C_0 X_1) \mathbf{Y} \\
&\quad + \varepsilon^2 (C_2 X_0 + C_1 X_1 + C_0 X_2) \mathbf{Y} + \cdots]^T \\
&\quad [B_0 + \varepsilon B_1 + \varepsilon^2 B_2 + \cdots] \\
&\quad [(\hat{\boldsymbol{\beta}} - \boldsymbol{\beta}_0) + \varepsilon (C_1 X_0 + C_0 X_1) \mathbf{Y} \\
&\quad + \varepsilon^2 (C_2 X_0 + C_1 X_1 + C_0 X_2) \mathbf{Y} + \cdots] \\
&= (\hat{\boldsymbol{\beta}} - \boldsymbol{\beta}_0)^T B_0 (\hat{\boldsymbol{\beta}} - \boldsymbol{\beta}_0) \\
&\quad + \varepsilon [(C_1 C_0^{-1} \hat{\boldsymbol{\beta}} + C_0 X_1 \mathbf{Y})^T B_0 (\hat{\boldsymbol{\beta}} - \boldsymbol{\beta}_0) \\
&\quad + (\hat{\boldsymbol{\beta}} - \boldsymbol{\beta}_0)^T B_1 (\hat{\boldsymbol{\beta}} - \boldsymbol{\beta}_0) \\
&\quad + (\hat{\boldsymbol{\beta}} - \boldsymbol{\beta}_0)^T B_0 (C_1 C_0^{-1} \hat{\boldsymbol{\beta}} + C_0 X_1 \mathbf{Y})] \\
&\quad + \varepsilon^2 [2 (C_2 C_0^{-1} \hat{\boldsymbol{\beta}} + C_1 X_1 \mathbf{Y} + C_0 X_2 \mathbf{Y})^T B_0 (\hat{\boldsymbol{\beta}} - \boldsymbol{\beta}_0) \\
&\quad + (\hat{\boldsymbol{\beta}} - \boldsymbol{\beta}_0)^T B_2 (\hat{\boldsymbol{\beta}} - \boldsymbol{\beta}_0) \\
&\quad + 2 (C_1 C_0^{-1} \hat{\boldsymbol{\beta}} + C_0 X_1 \mathbf{Y})^T B_1 (\hat{\boldsymbol{\beta}} - \boldsymbol{\beta}_0) \\
&\quad + (C_2 C_0^{-1} \hat{\boldsymbol{\beta}} + C_1 X_1 \mathbf{Y} B_0 (C_2 C_0^{-1} \hat{\boldsymbol{\beta}} C_1 X_1 \mathbf{Y} + C_0 X_2 \mathbf{Y} + C_0 X_2 \mathbf{Y})] \\
&\quad + \cdots \\
&= (\hat{\boldsymbol{\beta}} - \boldsymbol{\beta}_0)^T B_0 (\hat{\boldsymbol{\beta}} - \boldsymbol{\beta}_0) \\
&\quad + \varepsilon [2 \hat{\boldsymbol{\beta}}^T C_0^{-1} B_0 (\hat{\boldsymbol{\beta}} - \boldsymbol{\beta}_0) + 2 \mathbf{Y}^T X_1^T (\hat{\boldsymbol{\beta}} - \boldsymbol{\beta}_0) \\
&\quad + (\hat{\boldsymbol{\beta}} - \boldsymbol{\beta}_0)^T B_1 (\hat{\boldsymbol{\beta}} - \boldsymbol{\beta}_0)] \\
&\quad + \cdots
\end{aligned}
$$

Hence, up to the first order approximation, we can express $F(\varepsilon)$ as $F(\varepsilon) = N(\varepsilon)/D(\varepsilon)$ where

$$
\begin{aligned}
N(\varepsilon) &= \{(\hat{\boldsymbol{\beta}} - \boldsymbol{\beta}_0)^T B_0 (\hat{\boldsymbol{\beta}} - \boldsymbol{\beta}_0) + \varepsilon [2 \hat{\boldsymbol{\beta}}^T C_0^{-1} B_0 (\hat{\boldsymbol{\beta}} - \boldsymbol{\beta}_0) + 2 \mathbf{Y}^T X_1^T (\hat{\boldsymbol{\beta}} - \boldsymbol{\beta}_0) \\
&\quad + (\hat{\boldsymbol{\beta}} - \boldsymbol{\beta}_0)^T B_1 (\hat{\boldsymbol{\beta}} - \boldsymbol{\beta}_0)]\}/m
\end{aligned}
$$

and

$$
D(\varepsilon) = \frac{[SSE - \varepsilon \mathbf{Y}^T (X_1^T C_0 X_0 + X_0^T C_1 X_0 + X_0^T C_0 X_1) \mathbf{Y}]}{(n - m)}.
$$

The expression $F(\varepsilon)$ can be further expanded as

$$
\begin{aligned}
(3.14) \quad F(\varepsilon) &= \frac{(n-m)}{m}[(\hat{\boldsymbol{\beta}}-\boldsymbol{\beta}_0)^T B_0(\hat{\boldsymbol{\beta}}-\boldsymbol{\beta}_0) \\
&+ \varepsilon\{2\hat{\boldsymbol{\beta}}^T C_0^{-1} B_0(\hat{\boldsymbol{\beta}}-\boldsymbol{\beta}_0) \\
&+ 2\mathbf{Y}^T X_1^T(\hat{\boldsymbol{\beta}}-\boldsymbol{\beta}_0) \\
&+ (\hat{\boldsymbol{\beta}}-\boldsymbol{\beta}_0)^T B_1(\hat{\boldsymbol{\beta}}-\boldsymbol{\beta}_0)\}]\frac{1}{SSE} \\
&\quad [1+\frac{\varepsilon \mathbf{Y}^T(X_1^T C_0 X_0 + X_0^T C_1 X_0 + X_0^T C_0 X_1)\mathbf{Y}}{SSE}]+\cdots \\
&= F_0 + \varepsilon[(\frac{n-m}{m})\cdot 2 \cdot \frac{(\hat{\boldsymbol{\beta}}^T B_0 C_1 B_0(\hat{\boldsymbol{\beta}}-\boldsymbol{\beta}_0) + \mathbf{Y}^T X_1^T(\hat{\boldsymbol{\beta}}-\boldsymbol{\beta}_0)}{SSE} \\
&+ (\frac{n-m}{m})\cdot \frac{(\hat{\boldsymbol{\beta}}-\boldsymbol{\beta}_0)^T B_0(\hat{\boldsymbol{\beta}}-\boldsymbol{\beta}_0)}{SSE} \\
&\cdot \frac{\mathbf{Y}^T\{X_1^T C_0 X_0 + X_0^T C_1 X_0 + X_0^T C_0 X_1\}\mathbf{Y}}{SSE}]+\cdots
\end{aligned}
$$

where $F_0$, the usual $F$ statistic for testing $H_0$ versus $H_1$, is defined as

$$
(3.15) \quad F_0 = \frac{(\hat{\boldsymbol{\beta}}-\boldsymbol{\beta}_0)^T B_0(\hat{\boldsymbol{\beta}}-\boldsymbol{\beta}_0)/m}{SSE/(n-m)}.
$$

Lastly, we discuss the nature of the confidence set for $\boldsymbol{\beta}$ under the perturbed linear model (3.1). Obviously, under the assumption of normality and independence of errors, the $(1-\alpha)$ level confidence set for $\boldsymbol{\beta}$ is given by the following.

$$
\begin{aligned}
(3.16) \quad \mathcal{C}(\varepsilon) &= \{\boldsymbol{\beta} : (\hat{\boldsymbol{\beta}}(\varepsilon)-\boldsymbol{\beta})^T (X(\varepsilon)X^T(\varepsilon))(\hat{\boldsymbol{\beta}}(\varepsilon)-\boldsymbol{\beta}) \\
&\leq \frac{m}{n-m}\cdot F_{\alpha\,;\,m,n-m}\cdot SSE(\varepsilon)\} \\
&\equiv \{\boldsymbol{\beta} : F(\varepsilon) \leq F_{\alpha\,;\,m,n-m}\}.
\end{aligned}
$$

The dependence of $\mathcal{C}(\varepsilon)$ on $\varepsilon$ can then be studied based on the expansion of $F(\varepsilon)$ given in (3.14), and we have up to the first order

$$
\begin{aligned}
(3.17) \quad \mathcal{C}(\varepsilon) = \{\boldsymbol{\beta} : &\frac{(\hat{\boldsymbol{\beta}}-\boldsymbol{\beta})' B_0(\hat{\boldsymbol{\beta}}-\boldsymbol{\beta})/m}{SSE/(n-m)} \\
&+\varepsilon\cdot\frac{n-m}{n}[\frac{2\hat{\boldsymbol{\beta}}^T B_0 C_1 B_0(\hat{\boldsymbol{\beta}}-\boldsymbol{\beta})+\mathbf{Y}^T X_1^T(\hat{\boldsymbol{\beta}}-\boldsymbol{\beta})}{SSE} \\
&+ \frac{(\hat{\boldsymbol{\beta}}-\boldsymbol{\beta})^T B_0(\hat{\boldsymbol{\beta}}-\boldsymbol{\beta})}{SSE}\cdot\frac{\mathbf{Y}^T(X_1^T C_0 X_0 + X_0^T C_1 X_0 + X_0^T C_0 X_1)\mathbf{Y}}{SSE}] \\
&\leq F_{\alpha\,;\,m,n-m}\}.
\end{aligned}
$$

Case 2: $B_0 = X_0 X_0^T$ is singular. We now turn our attention to the case when $B_0$ is singular which is the singular perturbation situation and study the nature of dependence of the above estimates and test statistics on $\varepsilon$. We define $B(\varepsilon) = X(\varepsilon)X^T(\varepsilon)$, and consider first the case when

$$
(3.18) \quad B^{-1}(\varepsilon)X(\varepsilon) = \sum_{k=0}^{\infty}\varepsilon^k C_k^*.
$$

The above condition will not always hold, and necessary and sufficient conditions on the basic design matrices $X_0, X_1, \ldots$ can be developed under which such a representation holds by equating $X(\varepsilon)$ to $B(\varepsilon) \times \sum_{k=0}^{\infty} \varepsilon^k C_k^*$. For example, equating the first two terms yields conditions on the existence of $C_0^*$ and $C_1^*$ satisfying: $X_0 = (X_0 X_0^T) C_0^*$ and $X_1 = (X_0 X_0^T) C_1^* + (X_1 X_0^T + X_0 X_1^T) C_0^*$.

Under the above representation (3.18), let $\tilde{\boldsymbol{\beta}} = C_0^* \mathbf{Y}$ which is the limiting estimate of $\boldsymbol{\beta}$ as $\varepsilon$ approaches 0. It should be noted that $\tilde{\boldsymbol{\beta}} = B_0^G X_0 \mathbf{Y}$ for some generalized inverse $B_0^G$ of $B_0 = X_0 X_0^T$. To see this, note that under the representation (3.18),

$$X(\epsilon) = B(\epsilon) B^{-1}(\epsilon) X(\epsilon) = X(\epsilon) X^T(\epsilon) (\sum_{k=0}^{\infty} \epsilon^k C_k^*).$$

Taking the limit as $\epsilon \to 0$, we get

$$X_0 = X_0 X_0^T C_0^* = B_0 C_0^*.$$

Thus
$$C_0^* = B_0^G X_0$$

for some generalized inverse $B_0^G$ of $B_0$.

Returning to the error sum of squares $SSE(\varepsilon)$, recall that

$$SSE(\varepsilon) = \mathbf{Y}^T [I - X^T(\varepsilon) B^{-1}(\varepsilon) X(\varepsilon)] \mathbf{Y}$$

and

$$(3.19) \qquad \begin{aligned} P(\varepsilon) &= I_n - X^T(\varepsilon) B^{-1}(\varepsilon) X(\varepsilon) \\ &= \sum_{k=0}^{\infty} \varepsilon^k D_k^* \quad \text{(say)}, \end{aligned}$$

where we have used the representation (3.18). We define $\lim_{\varepsilon \to 0} SSE(\varepsilon) = \mathbf{Y}^T D_0^* \mathbf{Y}$, where

$$D_0^* = P(0) = I_n - X_0^T C_0^* = I_n - X_0^T B_0^G X_0,$$

which is always idempotent. Thus $SSE = SSE(\varepsilon = 0) = \mathbf{Y}^T [I_n - X_0 B_0^G X_0] \mathbf{Y}$ and because of the idempotency of $D_0^*$ we conclude that $\frac{SSE}{\sigma^2}$ follows a non-central chisquare distribution with non-centrality parameter $\beta' X(\epsilon) D_0^* X^T(\epsilon) \beta / \sigma^2$ and degrees of freedom $\nu$ given by

$$\begin{aligned} \nu &= \operatorname{rank}\{P(\varepsilon = 0)\} \\ &= n - \operatorname{rank}(X_0 B_0^G X_0) \\ &= n - r, \ r = \operatorname{rank}(X_0). \end{aligned}$$

Returning to the testing problem, let us define

$$F(\varepsilon) = \frac{(\tilde{\boldsymbol{\beta}}(\varepsilon) - \boldsymbol{\beta}_0)^T [X(\varepsilon) X^T(\varepsilon)][\tilde{\boldsymbol{\beta}}(\varepsilon) - \boldsymbol{\beta}_0]/m}{SSE(\varepsilon)/(n-m)}$$

$$(3.20) \qquad \xrightarrow[\varepsilon \to 0]{} \tilde{F} = \frac{(\tilde{\boldsymbol{\beta}} - \boldsymbol{\beta}_0)^T B_0 (\tilde{\boldsymbol{\beta}} - \boldsymbol{\beta}_0)/m}{SSE/(n-m)}.$$

We note that when $\varepsilon = 0$, for testing $H_0$: $X\boldsymbol{\beta} = X\boldsymbol{\beta}_0$, the F statistic, say $F_0$, is given by
$$F_0 = \frac{(\hat{\boldsymbol{\beta}} - \boldsymbol{\beta}_0)^T B_0 (\hat{\boldsymbol{\beta}} - \boldsymbol{\beta})/r}{SSE/(n-r)}$$
where
$$\hat{\boldsymbol{\beta}} = B_0^G X_0 \mathbf{Y}$$
and $r$ is the rank of $X_0$. Thus $\tilde{F}$ is a scalar multiple of $F_0$.

When the representation (3.18) assumed earlier under Case 2 does not hold, we would actually have a Laurent series expansion of $B^{-1}(\varepsilon)$ (as opposed to a Maclaurin series), resulting in

(3.21) $$B^{-1}(\varepsilon)X(\varepsilon) = \sum_{k=-s}^{\infty} \varepsilon^k C_k^*.$$

In this case
$$\tilde{\boldsymbol{\beta}}(\varepsilon) = B^{-1}(\varepsilon)X(\varepsilon)\mathbf{Y}$$
has no limit as $\varepsilon \downarrow 0$, but
$$E[\tilde{\boldsymbol{\beta}}(\varepsilon)] = B^{-1}(\varepsilon)X(\varepsilon)X^T(\varepsilon)\boldsymbol{\beta} = \boldsymbol{\beta}, \ \forall \varepsilon$$
$$\Rightarrow \quad \tilde{\boldsymbol{\beta}}(\varepsilon) \text{ is an unbiased estimate of } \boldsymbol{\beta}.$$

$$\begin{aligned}\text{Var}[\tilde{\boldsymbol{\beta}}(\varepsilon)] &= \sigma^2 B^{-1}(\varepsilon)X(\varepsilon)X^T(\varepsilon)B^{-1}(\varepsilon) \\ &= \sigma^2 B^{-1}(\varepsilon).\end{aligned}$$

Moreover, $\hat{\boldsymbol{\beta}}(\varepsilon)$ is also unbiased as $\varepsilon \to 0$ and
$$\begin{aligned}\text{Var}[\hat{\boldsymbol{\beta}}(\varepsilon)] &\sim \sigma^2 B_0^G X_0 X_0^T B_0^G \\ &= \sigma^2 B_0^G.\end{aligned}$$

Obviously, $\tilde{\boldsymbol{\beta}}(\varepsilon)$ and $\hat{\boldsymbol{\beta}}$ are equivalent for small $\varepsilon$ if $B^{-1}(\varepsilon) \xrightarrow[\varepsilon \downarrow 0]{} B_0^G$.

**Remark 3.1.** Interestingly enough, as proved in Avrachenkov [2], the projection matrix $P(\varepsilon)$ is uniformly bounded in $\varepsilon > 0$, and it does admit a Maclaurin series expansion at $\varepsilon = 0$, irrespective of the nature of $B^{-1}(\varepsilon)X(\varepsilon)$. Incidentally, the following simple example demonstrates that the projection matrix $P(\varepsilon)$ can be even independent of $\varepsilon$.

**Example 3.1.** Take
$$X_0^T = \begin{pmatrix} 2 & 0 \\ 2 & 0 \\ 1 & 0 \\ 1 & 0 \end{pmatrix}, \quad X_1^T = \begin{pmatrix} 0 & 1 \\ 0 & 1 \\ 0 & 1 \\ 0 & 1 \end{pmatrix}.$$

Then
$$X^T(\epsilon) = \begin{pmatrix} 2 & \epsilon \\ 2 & \epsilon \\ 1 & \epsilon \\ 1 & \epsilon \end{pmatrix}.$$

One might imagine that the two columns of $X^T(\epsilon)$ are generated by two sensors: the first one is accurate but the second one systematically overestimates the true reading of 0 by the amount $\epsilon$. If

$$\tilde{X}^T = \begin{pmatrix} 2 & 1 \\ 2 & 1 \\ 1 & 1 \\ 1 & 1 \end{pmatrix},$$

then we have

$$I_4 - P(\epsilon) = \tilde{X}^T(\tilde{X}\tilde{X}^T)^{-1}\tilde{X}.$$

Thus $P(\epsilon)$ is free of $\epsilon$. Hence, in this example, the systematic bias leaves the projection matrix (and hence also the residual sum of squares) unchanged.

**Remark 3.2.** Incidentally, an estimate of $\varepsilon$ can be provided by minimizing the error sum of squares, $SSE(\varepsilon)$, which is essentially the maximum likelihood estimate of $\varepsilon$. Recalling from (3.11) the representation of $SSE(\varepsilon)$ and keeping terms up to $\varepsilon^2$, we can verify that $SSE(\varepsilon)$ is convex in $\varepsilon$, and hence the minimizer solution is given by

(3.22)
$$\hat{\varepsilon} = \frac{\mathbf{Y}^T[X_1^TC_0X_0 + X_0^TC_1X_0 + X_0^TC_0X_1]\mathbf{Y}}{2\mathbf{Y}^T[X_2^TC_0X_0 + X_0^TC_2X_0 + X_0^TC_0X_2 + X_1^TC_1X_0 + X_0^TC_1X_1 + X_1^TC_0X_1]\mathbf{Y}}.$$

When the perturbation parameter $\varepsilon$ is unknown, the point estimate $\hat{\varepsilon}$ provides information about its magnitude. It is quite likely that the parameters of primary interest in a perturbed linear model are $\boldsymbol{\beta}$ and $\sigma^2$. After obtaining $\hat{\varepsilon}$, the maximum likelihood estimate of $\boldsymbol{\beta}$ can be computed as $\hat{\boldsymbol{\beta}}(\hat{\varepsilon})$. Similarly $\sigma^2$ can be estimated.

## 4. Perturbation in principal component analysis

In this section we discuss some effects of perturbation in principal component analysis (PCA). It is well known that in the standard PCA approach, one tries to reduce the dimension of a random $p \times 1$ vector $\mathbf{y}$ with a mean vector $\mu$ and a dispersion matrix $\Sigma$ by successively taking suitable orthogonal linear combinations of $\mathbf{y}$ such that the resultant linear combinations, known as principal components, have decreasing variances; see, for example, Anderson [1].

Generalizing the model (3.1) to a multivariate set up, we write

(4.1) $$\mathbf{Y}_{n\times p} = \mathbf{X}^T(\varepsilon)_{n\times m}\beta_{m\times p} + \psi_{n\times p}$$

where $\mathbf{Y}$ is a matrix of observations, the matrix $\mathbf{X}(\varepsilon)$ is assumed to be of full column rank, and the rows of $\psi$ are assumed to be independently distributed with mean vector $\mathbf{0}$ and a common dispersion matrix $\Sigma$. An estimate of the dispersion matrix $\Sigma$ is then obtained as

(4.2) $$\hat{\Sigma}(\varepsilon) = \mathbf{Y}^T P(\varepsilon)\mathbf{Y}$$

where, as before, $P(\varepsilon)$ is the projection matrix defined in (3.5). The sample principal components $\mathbf{d}_1^T\mathbf{y}$, $\mathbf{d}_2^T\mathbf{y}$, ... are then obtained by solving the equations:

(4.3) $$\hat{\Sigma}(\varepsilon)\mathbf{d}(\varepsilon) = \lambda(\varepsilon)\mathbf{d}(\varepsilon)$$

where $\lambda(\varepsilon)$'s are the eigenvalues of $\hat{\Sigma}(\varepsilon)$. When $\varepsilon = 0$, which means absence of any perturbation, the usual PCA applies. When $\varepsilon > 0$, the sample eigenvalues and the eigenvectors will depend on the perturbation parameter $\varepsilon$. There are quite a few papers dealing with the behaviors of such eigenvalues and eigenvectors as functions of $\varepsilon$ from a deterministic point of view.

In our context, the statistical properties of the random eigenvalues and eigenvectors can be studied based on such behaviors. It turns out that while the behaviors of the eigenvectors are smooth, those of the eigenvalues are not!

Following Remark 2.4, we can write $P(\varepsilon)$ as

$$P(\varepsilon) = P_0 + \varepsilon P_1 + \varepsilon^2 P_2 + \cdots \tag{4.4}$$

which readily yields

$$\hat{\Sigma}(\varepsilon) = S_0 + \varepsilon S_1 + \varepsilon^2 S_2 + \cdots \tag{4.5}$$

where $S_0 = SSE$. The eigenvalues $\lambda(\varepsilon)$ then satisfy the determinantal equation:

$$|S_0 + \varepsilon S_1 + \varepsilon^2 S_2 + \cdots - \lambda(\varepsilon) I_p| = 0 \tag{4.6}$$

and once the $\lambda(\varepsilon)$'s have been determined, the eigen vectors $\mathbf{d}(\varepsilon)$'s are obtained from the equations given in (4.3).

In the special case when $P(\varepsilon) = P_0 + \varepsilon P_1$, the equation (4.6) reduces to

$$|S_0 + \varepsilon S_1 - \lambda(\varepsilon) I_p| = 0 \tag{4.7}$$

which when further specialized to $p = 2$ yields

$$\lambda_1(\varepsilon) + \lambda_2(\varepsilon) = tr(S_0) + \varepsilon tr(S_1) \quad \lambda_1(\varepsilon)\lambda_2(\varepsilon) = |S_0 + \varepsilon S_1|.$$

Since $(\lambda_1(\varepsilon) - \lambda_2(\varepsilon))^2 = (\lambda_1(\varepsilon) + \lambda_2(\varepsilon))^2 - 4\lambda_1(\varepsilon)\lambda_2(\varepsilon)$, up to the first order of approximation, we get

$$\lambda_1(\varepsilon) - \lambda_2(\varepsilon) \sim A + \varepsilon \frac{[(s_{111} - s_{122})(s_{011} - s_{022}) + 4s_{012}s_{112}]}{A} \tag{4.8}$$

where

$$A = [(s_{011} - s_{022})^2 + 4s_{012}^2]^{1/2}. \tag{4.9}$$

In the above, we have used the notation

$$S_0 = \begin{bmatrix} s_{011} & s_{012} \\ s_{012} & s_{022} \end{bmatrix}$$

$$S_1 = \begin{bmatrix} s_{111} & s_{112} \\ s_{112} & s_{122} \end{bmatrix}.$$

The equations (4.8) and (4.9) can be simultaneously solved to get first order approximations of $\lambda_1(\varepsilon)$ and $\lambda_2(\varepsilon)$. Once this is done, first order approximations of the two eigenvectors $\mathbf{d}_1(\varepsilon) \sim \mathbf{d}_{10} + \varepsilon \mathbf{d}_{11}$ and $\mathbf{d}_2(\varepsilon) \sim \mathbf{d}_{20} + \varepsilon \mathbf{d}_{21}$ are obtained by solving the equations:

$$[S_0 + \varepsilon S_1](\mathbf{d}_{10} + \varepsilon \mathbf{d}_{11}) = \lambda_1(\varepsilon)(\mathbf{d}_{10} + \varepsilon \mathbf{d}_{11}) \tag{4.10}$$

(4.11) $$[S_0 + \varepsilon S_1](\mathbf{d}_{20} + \varepsilon \mathbf{d}_{21}) = \lambda_2(\varepsilon)(\mathbf{d}_{20} + \varepsilon \mathbf{d}_{21})$$

Writing $\lambda_1(\varepsilon) \sim \lambda_{10} + \varepsilon \lambda_{11}$, we readily get $\mathbf{d}_{10}$ and $\mathbf{d}_{11}$ as solutions of the equations:

$$S_0 \mathbf{d}_{10} = \lambda_{10} \mathbf{d}_{10} S_0 \mathbf{d}_{11} + S_1 \mathbf{d}_{10} = \lambda_{11} \mathbf{d}_{10} + \lambda_{10} \mathbf{d}_{11}.$$

Likewise, the components of the other eigenvector $\mathbf{d}_2(\varepsilon)$ can be obtained.

It would indeed be challenging to study the statistical properties of the resultant eigenvectors and eigenvalues.

## 5. Perturbation in factor analysis

In this section we discuss some consequences of perturbation in the context of factor analysis (FA).

Let $\mathbf{Y} = (Y_1, Y_2, \ldots, Y_p)$ be a vector of *manifest* variables. Then the FA model postulates that (Christensen [7] and Johnson and Wichern [10]) the manifest variables are linear functions of some random latent variables plus a residual term. The functional relationship is typically expressed as

(5.1)
$$\begin{aligned}
Y_1 &= \lambda_{11} f_1 + \lambda_{12} f_2 + \cdots + \lambda_{1k} f_k + \psi_1 \\
Y_2 &= \lambda_{21} f_1 + \lambda_{22} f_2 + \cdots + \lambda_{2k} f_k + \psi_2 \\
&\ldots \\
Y_p &= \lambda_{p1} f_1 + \lambda_{p2} f_2 + \cdots + \lambda_{pk} f_k + \psi_p
\end{aligned}$$

where $f_1, f_2, \ldots, f_k$ represent the random latent or common factors and $(\psi_1, \psi_2, \ldots, \psi_p)$ denote the residual or error terms. Here the elements $\lambda$'s are known as factor loadings, the latent variables are assumed to have a normal distribution with mean vector 0 and a dispersion matrix $\Phi$, and the residuals are assumed to be independent of $f$'s, and independently normally distributed with mean 0 and a variance-covariance matrix $\Psi = \mathrm{diag}(\psi_{11}, \psi_{22}, \ldots, \psi_{pp})$. The dispersion matrix $\Sigma$ of $\mathbf{Y}$ is then given by

(5.2) $$\Sigma = \Gamma \Phi \Gamma^T + \Psi.$$

The literature on FA usually assumes that the latent variables are orthogonal, resulting in $\Phi = I_k$ and hence the simplification $\Sigma = \Gamma \Gamma^T + \Psi$.

Since $E\mathbf{Y} = \mathbf{0}$, the sample covariance matrix $S$ can be used to estimate $\Sigma$, and the parameters in $\Gamma$ and $\Psi$ are estimated by solving

(5.3) $$S = \hat{\Gamma}\hat{\Gamma}^T + \hat{\Psi},$$

subject to the condition that $\hat{\Gamma}^T(\hat{\Psi})^{-1}\hat{\Gamma}$ is diagonal. Iterative methods can be used to solve the above equations.

The method of maximum likelihood, on the other hand, is based on maximizing the likelihood or its logarithm given by

(5.4) $$\mathcal{L}(\Gamma, \Psi | S) \sim -[\ln\{\det(\Sigma)\} + \mathrm{tr}(S\Sigma^{-1})]$$

which is directly a function of factor loadings and the error variances in view of (5.2). An old but most successful algorithm due to Joreskog [11] which runs in two steps, first by maximizing $\mathcal{L}$ with respect to $\Gamma$ for a fixed $\Psi$, and then by maximizing with respect to $\Psi$, can be used to derive the maximum likelihood estimates. Details are omitted.

Returning to the perturbation formulation, we can introduce the perturbation parameter $\varepsilon$ in the form

$$(5.5) \qquad \Phi(\varepsilon) = I_k + \varepsilon\Phi_1 + \varepsilon^2\Phi_2 + \cdots$$

which takes us away from the usual assumption of independence of the latent variables and emphasizes that there may be some uncertainties in this assumption. We can then rewrite $\Sigma$ as

$$(5.6) \qquad \begin{aligned} \Sigma(\varepsilon) &= \Gamma\Phi(\varepsilon)\Gamma^T + \Psi \\ &= \Gamma[I_k + \varepsilon\Phi_1 + \varepsilon^2\Phi_2 + \cdots]\Gamma^T + \Psi \\ &= [\Gamma\Gamma^T + \Psi] + \varepsilon\Gamma\Phi_1\Gamma^T + \varepsilon^2\Gamma\Phi_2\Gamma^T + \cdots \end{aligned}$$

which can be used to compute and expand $\Sigma(\varepsilon)^{-1}$ and also $\ln|\Sigma(\varepsilon)|$ in powers of $\varepsilon$. The maximum likelihood estimates of the factor loadings (elements of $\Gamma$) and the variances (elements of $\Psi$) would then naturally depend on the perturbation parameter $\varepsilon$, and it is indeed possible to study the effect of $\varepsilon$ on these estimates along the same lines as in the previous sections. Details are omitted.

## 6. Applications

In this section we provide an application of the theory developed in this paper. Consider the nonlinear regression of $Y$ on $\mathbf{X}$ given by (Gallant [8], Chapter 1, Example 1)

$$(6.1) \qquad Y_i = \theta_1 x_{1i} + \theta_2 x_{2i} + \theta_0 e^{\epsilon x_{3i}} + e_i, \ \ i = 1, \ldots, n$$

where the $Y_i$'s are the responses of a *treatment-control* design, $x_1 = 1$ and 0 represent, respectively, treatment and control scenarios, $x_2$ represents a variable with values between 1 and 2, and $x_3$ is the variable showing age of the experimental material. The errors $e_i$'s are assumed to be (normally) distributed with mean 0 and variance $\sigma^2$. For $\epsilon = 0$, this of course is a simple linear regression of $Y$ on $(x_1, x_2)$ for which standard inference applies. For $\epsilon > 0$, this is an application of a non-linear regression set up for which drawing appropriate inference about the parameters is rather complex. An expansion of $e^{\epsilon x_{3i}}$ readily yields the matrices $X_0, X_1, X_2, \ldots$ given in (3.2), from which the relevant quantities $\hat{\beta}(\epsilon)$, $SSE(\epsilon)$, $F(\epsilon)$ and $\mathcal{C}(\epsilon)$ can be computed up to any order of powers of $\epsilon$.

The data given in Gallant [8], page 4, has all the $x_{2i}$ values equal to one. When this is the case, it is clear from (6.1) that $\theta_2$ and $\theta_0$ are not identifiable when $\epsilon = 0$. Consequently, in order to illustrate our methodology, we shall avoid the choice $x_{2i} = 1$ for all $i$. In the data given in Table 1, the $x_{2i}$-values were randomly chosen from a uniform distribution on $(1, 2)$. Also given in Table 1 are four data sets, simulated based on the model (6.1). The data sets share the same set of values of the covariates $x_1$, $x_2$ and $x_3$, and differ in the values of the response variable $Y$. In each case, $n = 30$.

To see the effect of $\epsilon$ on the standard analysis, we note that for the given data sets, $X_o^T : 30 \times 3$ is a matrix consisting of the first three columns from Table 1, $X_1^T : 30 \times 3$ is a matrix whose first column are the elements $x_3$ and the remaining two columns are null vectors, and lastly, $X_2^T : 30 \times 3$ is a matrix whose first column elements are $(0.5)$ times the squares of $x_3$ and the remaining two columns are null vectors. Using these basic matrices, we readily compute

$$B_0 = X_0 X_0^T = \begin{pmatrix} 30 & 15 & 44.8534 \\ 15 & 15 & 21.7371 \\ 44.8534 & 21.7371 & 69.3119 \end{pmatrix}.$$

$$C_0 = B_0^{-1} = \begin{pmatrix} 1.1523 & -0.1314 & -0.7045 \\ -0.1314 & 0.1372 & 0.0420 \\ -0.7045 & 0.0420 & 0.4571 \end{pmatrix}.$$

$$B_1 = X_1 X_0^T + X_0 X_1^T = \begin{pmatrix} 294.62 & 74.55 & 214.3236 \\ 74.55 & 0 & 0 \\ 214.3236 & 0 & 0 \end{pmatrix}.$$

$$C_1 = -C_0 B_1 C_0 = \begin{pmatrix} -20.6598 & 1.3236 & 9.3912 \\ 1.3236 & -0.0332 & -0.4397 \\ 9.3912 & -0.4397 & -3.7610 \end{pmatrix}.$$

$$B_2 = X_2 X_0^T + X_1 X_1^T + X_0 X_2^T = \begin{pmatrix} 2035.2007 & 264.923 & 738.250 \\ 264.923 & 0 & 0 \\ 738.250 & 0 & 0 \end{pmatrix}.$$

$$C_2 = -B_0^{-1}(B_2C_0 + B_1C_1) = \begin{pmatrix} -164.4712 & 32.577 & 147.827 \\ 32.577 & -4.3659 & -22.4049 \\ 147.8266 & -22.4049 & -110.1712 \end{pmatrix}.$$

TABLE 1
*Data sets*

| $x_0$ | $x_1$ | $x_2$ | $x_3$ | $Y_0$ | $Y_1$ | $Y_2$ | $Y_3$ |
|---|---|---|---|---|---|---|---|
| 1 | 1 | 1.3420 | 6.28 | 6.92 | 7.0451 | 7.6532 | 8.6632 |
| 1 | 0 | 1.5813 | 9.86 | 4.21 | 4.4267 | 5.4939 | 7.5804 |
| 1 | 1 | 1.1043 | 9.11 | 4.34 | 4.5297 | 5.4928 | 7.3125 |
| 1 | 0 | 1.6867 | 8.43 | 3.81 | 3.9869 | 4.8595 | 6.4576 |
| 1 | 1 | 1.5164 | 8.11 | 3.99 | 4.1633 | 4.9944 | 6.4946 |
| 1 | 0 | 1.5672 | 1.82 | 3.95 | 3.9820 | 4.1358 | 4.3445 |
| 1 | 1 | 1.9644 | 6.58 | 6.28 | 6.4206 | 7.0637 | 8.1464 |
| 1 | 0 | 1.5411 | 5.02 | 5.63 | 5.7309 | 6.1986 | 6.9320 |
| 1 | 1 | 1.0064 | 6.52 | 4.05 | 4.1857 | 4.8217 | 5.8897 |
| 1 | 0 | 1.8726 | 3.75 | 5.53 | 5.6102 | 5.9462 | 6.4438 |
| 1 | 1 | 1.0314 | 9.86 | 5.09 | 5.3012 | 6.3684 | 8.4550 |
| 1 | 0 | 1.9190 | 7.31 | 6.16 | 6.3080 | 7.0388 | 8.3106 |
| 1 | 1 | 1.6507 | 0.47 | 5.80 | 5.8132 | 5.8514 | 5.9001 |
| 1 | 0 | 1.7083 | 0.07 | 4.98 | 4.9816 | 4.9873 | 4.9943 |
| 1 | 1 | 1.1261 | 4.07 | 4.97 | 5.0493 | 5.4176 | 5.9709 |
| 1 | 0 | 1.1693 | 4.61 | 5.38 | 5.4791 | 5.9032 | 6.5561 |
| 1 | 1 | 1.8063 | 0.17 | 7.19 | 7.1955 | 7.2092 | 7.2264 |
| 1 | 0 | 1.7086 | 6.99 | 5.19 | 5.3309 | 6.0228 | 7.2096 |
| 1 | 1 | 1.4324 | 4.39 | 6.20 | 6.2895 | 6.6907 | 7.3021 |
| 1 | 0 | 1.5265 | 0.39 | 5.14 | 5.1441 | 5.1757 | 5.2158 |
| 1 | 1 | 1.7009 | 4.73 | 4.37 | 4.4670 | 4.9038 | 5.5798 |
| 1 | 0 | 1.5807 | 9.42 | 3.82 | 4.0196 | 5.0253 | 6.9523 |
| 1 | 1 | 1.5538 | 8.90 | 5.38 | 5.5706 | 6.5055 | 8.2547 |
| 1 | 0 | 1.4150 | 3.02 | 3.12 | 3.1812 | 3.4459 | 3.8250 |
| 1 | 1 | 1.8566 | 0.77 | 5.06 | 5.0730 | 5.1360 | 5.2176 |
| 1 | 0 | 1.5010 | 3.31 | 4.08 | 4.1479 | 4.4406 | 4.8653 |
| 1 | 1 | 1.1584 | 4.51 | 5.72 | 5.8164 | 6.2300 | 6.8639 |
| 1 | 0 | 1.3310 | 2.65 | 3.89 | 3.9459 | 4.1756 | 4.4991 |
| 1 | 1 | 1.4981 | 0.08 | 6.12 | 6.1169 | 6.1234 | 6.1314 |
| 1 | 0 | 1.0008 | 6.11 | 2.84 | 2.9685 | 3.5571 | 4.5271 |

Hence, using (3.22), an estimate of $\epsilon$ for the four data sets is obtained as $\hat{\epsilon}_0 = 0.0505$, $\hat{\epsilon}_1 = 0.0195$, $\hat{\epsilon}_2 = -0.0227$, $\hat{\epsilon}_3 = -0.1861$. This results in the *estimated* regression coefficients and the *estimated F*-values (for testing the null hypothesis $H_0 : \theta_0 = \theta_1 = \theta_2 = 1$), up to the first order, as

*Estimated regression coefficients (standard errors) and F-values*

| Data Set | $\hat{\theta}_0$ | $\hat{\theta}_1$ | $\hat{\theta}_2$ | F |
|---|---|---|---|---|
| 1 | 2.0153 (0.9632) | 1.0642 (0.3329) | 1.6228 (0.6068) | 48.7884 |
| 2 | 2.2057 (0.9513) | 1.0629 (0.3287) | 1.5645 (0.5993) | 55.2427 |
| 3 | 3.0969 (0.9724) | 1.0545 (0.3360) | 1.2967 (0.6126) | 80.1393 |
| 4 | 4.6550 (1.3455) | 1.0537 (0.4650) | 0.8177 (0.8476) | 73.7195 |

On the other hand, taking $\epsilon = 0$, the *estimated* regression coefficients and the *F*-values, under the same null hypothesis $H_0 : \theta_0 = \theta_1 = \theta_2 = 1$ as above, are obtained in the four cases as

*Estimated regression coefficients (standard errors) and F-values when $\epsilon = 0$*

| Data Set | $\hat{\theta}_0$ | $\hat{\theta}_1$ | $\hat{\theta}_2$ | F |
|---|---|---|---|---|
| 1 | 2.0142 (0.9605) | 1.0640 (0.3327) | 1.6235 (0.6064) | 48.8080 |
| 2 | 2.2016 (0.9499) | 1.0618 (0.3283) | 1.5673 (0.5984) | 55.3330 |
| 3 | 3.0930 (0.9725) | 1.0549 (0.3362) | 1.2986 (0.6126) | 80.0064 |
| 4 | 4.6511 (1.3450) | 1.0540 (0.4650) | 0.8199 (0.8472) | 73.7090 |

It is interesting to note that there is practically no change in the estimated regression coefficients or the *F*-values when $\epsilon$ is zero or non-zero. Thus we conclude that for inference concerning the model (6.1), the non-linearity due to the presence of $\epsilon$ is inconsequential; we might as well assume that $\epsilon = 0$, and carry out the usual linear model analysis.

**Acknowledgments.** Our sincere thanks are due to Professor Edsel Peña for his constructive suggestions which led to improved presentation of results.

## References

[1] ANDERSON, T. W. (1984). *An Introduction to Multivariate Statistical Analysis*, 2nd ed. Wiley, New York. MR0771294
[2] AVRACHENKOV, K. E. (1999). Analytic perturbation theory and its applications. Doctoral dissertation, University of South Australia.
[3] AVRACHENKOV, K. E. AND HAVIV, M. (2003). Perturbation of null spaces with application to the eigenvalue problem and generalized inverses. *Linear Algebra Appl.* **369** 1–25. MR1988476
[4] AVRACHENKOV, K. E., HAVIV, M. AND HOWLETT, P. G. (2001). Inversion of analytic matrix functions that are singular at the origin. *SIAM J. Matrix Anal. Appl.* **22** 1175–1189. MR1824064
[5] CALAFIORE, G. AND GHAOUT, L. E. (2001). Robust maximum likelihood estimation in the linear model. *Automatica* **37** 573–580. MR1832529

[6] CHANDRASEKARAN, G. G., GU, M. AND SAYED, A. H. (1998). Parameter estimation in the presence of bounded data uncertainties. *SIAM J. Matrix Anal. Appl.* **19** 235–252. MR1609972

[7] CHRISTENSEN, R. (2001). *Advanced Linear Modeling*, 2nd ed. Springer, New York. MR1880721

[8] GALLANT, A. R. (1987). *Nonlinear Statistical Models*. Wiley, New York. MR0921029

[9] GHAOUI, L. E. AND LEBRET, H. (1997). Robust solutions to least squares problems with uncertain data. *SIAM J. Matrix Anal. Appl.* **18** 1035–1064. MR1472008

[10] JOHNSON, R. A. AND WICHERN, D. W. (2002). *Applied Multivariate Statistical Analysis*, 5th ed. Prentice Hall, New York.

[11] JORESKOG, K. G. (1967). Some contributions to maximum likelihood factor analysis. *Psychometrika* **32** 443–482. MR0221659

[12] KARIYA, T. AND SINHA, B. K. (1989). *Robustness of Statistical Tests*. Academic Press, New York. MR0996634

[13] KATO, T. (1980). *Perturbation Theory for Linear Operators*, 2nd ed. Springer, New York.